\documentclass[11pt]{article}
\normalsize
\usepackage{amsfonts}
\usepackage{amsmath}
\usepackage{amssymb}

\newtheorem{theorem}{Theorem}

\newtheorem{lemma}[theorem]{Lemma}

\title{ On the classification of binary self-dual $[44,22,8]$ codes
with an automorphism of order $3$ or $7$}
\author{Stefka Bouyuklieva,\\
Faculty of Mathematics and Informatics,\\ Veliko Tarnovo University, Bulgaria,\\
Nikolay Yankov and Radka Russeva,\\ Faculty of Mathematics and
Informatics,\\ Shumen University, Bulgaria\\}
\date{}

\begin{document}
\maketitle

\begin{abstract}
All binary self-dual $[44,22,8]$ codes with an automorphism of
order $3$ or $7$ are classified. In this way we complete the
classification of extremal self-dual codes of length 44 having an
automorphism of odd prime order.
\end{abstract}


\section{Introduction}

Let $\mathbb{F}_q^n$ be the $n$-dimensional vector space over the
field $\mathbb{F}_q$ of $q$ elements. A \emph{linear $[n,k]$ code}
$C$ is a $k$-dimensional subspace of $\mathbb{F}_q^n$. The elements
of $C$ are called \emph{codewords}. The \emph{weight} of a vector
$v\in\mathbb{F}_q^n$ (denoted by $\mbox{wt}(v)$) is the number of
its non-zero coordinates.
 The \emph{minimum weight} $d$ of $C$ is the
smallest weight among all nonzero weights of codewords of $C$; a
code $C$ with minimum weight $d$ is called an $[n,k,d]$ code. A
matrix whose rows form a basis of $C$ is called a \emph{generator
matrix} of this code. The weight enumerator $W(y)$ of a code $C$
is given by $W(y)=\sum_{i=0}^n A_iy^i$ where $A_i$ is the number
of codewords of weight $i$ in $C$. Two binary codes are called
\emph{equivalent} if one can be obtained from the other by a
permutation of coordinates. The permutation $\sigma\in S_n$ is an
\emph{automorphism} of $C$, if $C=\sigma(C)$ and the set of all
automorphisms of $C$ forms a group called the \emph{automorphism
group} of $C$, which is denoted by $Aut(C)$ in this paper.

Let $(u,v)\in\mathbb{F}_q$ for $u,v\in\mathbb{F}_q^n$ be an inner
product in $\mathbb{F}_q^n$. The \emph{dual code} of an $[n,k]$
code $C$ is $C^{\perp}=\{u \in \mathbb{F}_q^n \mid (u,v)=0$ for
all $v \in C \}$ and $C^{\perp}$ is a linear $[n,n-k]$ code. If $C
\subseteq C^{\perp}$, $C$ is termed \emph{self-orthogonal}, and if
$C = C^{\perp}$, $C$ is self-dual. We call a binary code
\emph{self-complementary} if it contains the all-ones vector.
Every binary self-dual code is self-complementary. If
$u=(u_1,\cdots, u_n)$, $v = (v_1,\cdots,v_n)\in\mathbb{F}_2^n$
then $(u,v)=\sum_{i=1}^{n}{u_iv_i}\in\mathbb{F}_2$. It was shown
in \cite{Rains98} that the minimum weight $d$ of a binary
self-dual code of length $n$ is bounded by $d \le 4[n/24]+4$,
unless $n \equiv 22 \pmod{24}$ when $d \le 4[n/24]+6$. We call a
self-dual code meeting this upper bound {\em extremal}.

 In this paper, we consider extremal binary self-dual $[44,22,8]$
codes. All the odd primes $p$ dividing the order of the
automorphism group of such a code are 11, 7, 5, and 3 \cite{160}.
The codes with automorphisms of order 11 and 5 are classified in
\cite{160}, \cite{164}, \cite{St42-44}, and \cite{Stef-SO}.
Unfortunately we noticed that there are some omissions in the
classification of the codes with automorphisms of order 7 given in
\cite{RN-SMB'08}. That's why we focus on the automorphisms of
orders 3 and 7, and we complete the classification of $[44,22,8]$
self-dual codes having an automorphism of odd prime order.

As in the case of binary self-dual $[42,21,8]$ codes with an
automorphism of order 3, there are five different possibilities
for the number of independent cycles in the decomposition of the
automorphism, namely 6, 8, 10, 12, and 14 \cite{BYR}. Codes with
automorphisms of order 3 with 6 and 14 independent 3-cycles are
considered but not classified in \cite{Stef-SO} and
\cite{Nikolay3-14}, respectively. In this paper, we give the
classification of all self-dual $[44,22,8]$ codes having an
automorphism of order 3 or 7.  To do that we apply the method for
 constructing binary self-dual codes via an automorphism of odd
 prime order developed in \cite{Huff} and \cite{Yo}. We give a short description of
this method  in Section 2. In Section 3 and Section 4 we classify
the extremal self-dual codes of length 44 with an automorphism of
order 3 and 7, respectively. In Section 5 we present the full classification of the self-dual $[44,22,8]$ codes having automorphisms of odd prime order, and offer some open problems.

The weight enumerators of the extremal self-dual codes of length
$44$ are known (see \cite{CS}):
$$W_{44,1}(y) =1 + (44 + 4\beta)y^8 + (976-8\beta)y^{10} +
(12289-20\beta)y^{12}+\dots$$ for $10\leq \beta\leq 122$  and
$$W_{44,2}(y) =1 + (44 + 4\beta)y^8 + (1232-8\beta)y^{10} +
(10241-20\beta)y^{12}+\dots$$ for $0\leq \beta\leq 154$.

Codes exist for $W_{44,1}$ when $\beta= 10,\dots, 68$, 70, 72, 74,
82, 86, 90, 122 and for $W_{44,2}$ when $\beta = 0, \dots , 56$,
58, \dots , 62, 64, 66, 68, 70, 72, 74, 76, 82, 86, 90, 104, 154
(see \cite{FFTAHuff}).

\section{Construction Method}

Let $C$ be a binary self-dual code of length $n=44$ with an
automorphism $\sigma $ of prime order $p\geq 3$ with exactly $c$
independent $p$-cycles and $f=44-cp$ fixed points in its
decomposition. We may assume that
$$\sigma=(1,2,\cdots,p)(p+1,p+2,\cdots,2p)\cdots(p(c-1)+1,p(c-1)+2,\cdots,pc),$$
and say that $\sigma$ is of type $p$-$(c,f)$.

Denote the cycles of $\sigma $ by $\Omega_1,\ldots, \Omega_c$, and
the fixed points by $\Omega_{c+1},\ldots, \Omega_{c+f}$. Let
$F_\sigma(C)=\{ v \in C\mid v \sigma=v \}$ and
$$E_\sigma(C)=\{ v\in C\mid\mbox{wt}(v\vert \Omega_i)\equiv 0\pmod
2, i=1,\cdots,c+f\},$$ where $v\vert \Omega_i $ is the restriction
of $v$ on $\Omega_i $.

\begin{theorem} {\rm \cite{Huff}} The self-dual code $C$ is a direct sum of the subcodes
$F_\sigma(C)$  and $E_\sigma(C)$. These subcodes have
dimensions $\frac{c+f}{2}$ and $\frac{c(p-1)}{2}$, respectively.
\end{theorem}

Thus each choice of the codes $F_\sigma(C)$ and $E_\sigma(C)$
determines a self-dual code $C$. So for a given length all
self-dual codes with an automorphism $\sigma$ can be obtained.

We have that $v\in F_\sigma (C)$ if and only if $v\in C$ and $v$
is constant on each cycle. Let $\pi :F_\sigma (C)\to
\mathbb{F}_{2}^{c+f} $ be the projection map where if $v\in
F_\sigma (C)$, $(\pi(v))_i =v_j $ for some $j\in\Omega_i$,
$i=1,2,\ldots,c+f $.

Denote by $E_\sigma(C)^{*}$ the code $E_\sigma(C)$ with the last
$f$ coordinates deleted. So $E_\sigma(C)^{*}$ is a self-orthogonal
binary code of length $pc$ and dimension $c(p-1)/2$. For $v\in
E_\sigma (C)^*$ we let $v\vert\Omega_i=(v_0,v_1,\cdots,v_{p-1})$
correspond to the polynomial $v_0+v_1 x+\cdots +v_{p-1}x^{p-1}$
from ${\cal P}$, where ${\cal P}$ is the set of even-weight
polynomials in ${\cal R}_p=\mathbb{F}_2[x]/\langle x^p-1\rangle$.
Thus we obtain the map $\varphi:E_\sigma(C)^{*}\to {\cal P}^c $.
${\cal P}$ is a cyclic code of length $p$ with generator
polynomial $x-1$.  It is known that $\varphi(E_{\sigma}(C)^{*})$
is a submodule of the ${\cal P}$-module ${\cal P}^c$
\cite{Huff,Yorus}.

\begin{theorem}{\rm \cite{Yorus}}\label{thm2}
A binary $[n,n/2]$ code $C$ with an automorphism $\sigma$ is
self-dual if and only if the following two conditions hold:

\begin{itemize}
\item[(i)]
$C_{\pi}=\pi(F_\sigma(C))$ is a binary self-dual code of length
$c+f$,
\item[(ii)] for every two vectors $u, v$ from $C_{\varphi}=\varphi(E_{\sigma}(C)^{*})$
we have
$$u_1(x)v_1(x^{-1})+ u_2(x)v_2(x^{-1})+ \cdots + u_c(x)v_c(x^{-1}) =
0.$$
\end{itemize}
\end{theorem}

Let $x^p-1=(x-1)h_1(x)\dots h_s(x)$, where $h_1,\dots,h_s$ are
irreducible binary polynomials. If $g_j(x)=(x^p-1)/h_j(x)$, and
$I_j=\langle g_j(x)\rangle$ is the ideal in ${\cal R}_p$,
generated by $g_j(x)$, then $I_j$ is a fields with
$2^{deg(h_j(x))}$ elements, $j=1, 2, \dots, s$, and ${\cal P} =
I_1\oplus I_2\oplus \cdots \oplus I_s$ \cite{Handbook}.

\begin{lemma}\label{M1M2M3} {\rm
\cite{Yorus}} Let $M_j=\{ u\in \varphi(E_\sigma(C)^{*})\vert
u_i\in I_j,i=1,2,\ldots,c\}$, $j=1, 2, \dots, s$. Then

1) $M_j$ is a linear space over $I_j$, $j=1, 2, \dots, s$;

2) $C_\varphi=\varphi(E_\sigma (C)^{*})=M_1\oplus M_2\oplus \cdots
\oplus M_s$ (direct sum of  ${\cal P}$-submodules);

3) If $C$ is a self-dual code, then $\sum\limits_{j =
1}^s{dim_{I_j}M_j = cs/2}$.
  \end{lemma}

In the case, when 2 is a primitive root modulo $p$, ${\cal P}$ is
a field with $2^{p-1}$ elements and the following theorem holds

\begin{theorem}{\rm \cite{Huff}}
\label{primroot} Let $2$ be a primitive root modulo $p$. Then the
binary code $C$ with an automorphism $\sigma$ is self-dual iff the
following two conditions hold:
\begin{itemize}
  \item[(i)]$C_\pi$ is a self-dual binary code of length $c + f$;
  \item[(ii)]$C_\varphi$ is a self-dual code of length $c$ over the field ${\cal P}$ under the inner product
$(u,v)=\sum_{i=1}^{c}{u_iv_i^{(p-1)/2}}.$
\end{itemize}
\end{theorem}

Let ${\cal B}$, respectively ${\cal D}$, be the largest subcode of
$C_\pi$ whose support is contained entirely in the left $c$,
respectively, right $f$, coordinates. Suppose ${\cal B}$ and
${\cal D}$ have dimensions $k_1$ and $k_2$, respectively. Let $k_3
= k - k_1 - k_2$. Then there exists a generator matrix for $C_\pi$
in the form
\begin{equation}\label{matGpi}
G_{\pi} =\left(\begin{array}{cc}
B & \ O\\
O & \ D\\
E & \ F
\end{array}\right)\end{equation}
where $B$ is a $k_1\times c$ matrix with $gen ({\cal B}) = [B \
O]$, $D$ is a $k_2\times f$ matrix with $gen ({\cal D}) = [O \
D]$, $O$ is the appropriate size zero matrix, and $[E \ F]$ is a
$k_3\times n$ matrix.
 Let ${\cal B}^*$ be the code of
length $c$ generated by $B$, ${\cal B}_E$ the code of length $c$
generated by the rows of $B$ and $E$, ${\cal D}^*$  the code of
length $f$ generated by $D$, and ${\cal D}_F$ the code of length
$f$ generated by the rows of $D$ and $F$. The following theorem is
 a modification of Theorem 2 from \cite{Handbook-Pless}.

\begin{theorem}\label{structure}
 With the notations of the previous paragraph
\begin{itemize}
  \item[(i)]$k_3 = \mbox{rank}(E) = \mbox{rank}(F)$,
  \item[(ii)]$k_2 = k + k_1 - c=k_1+\frac{f-c}{2}$, and
  \item[(iii)]${\cal B}_E^{\perp} = {\cal B}^*$ and ${\cal D}_F^{\perp} = {\cal D}^*$.
\end{itemize}
\end{theorem}

\section{Extremal Self-Dual Codes of Length 44 with an Automorphism of
Order 3}

Using Theorem \ref{primroot}, as 2 is a primitive root modulo 3,
${\cal P}$ is a field with 4 elements. We have that ${\cal P}=\{0,
e=x+x^2, w=1+x^2, w^2=1+x\}\cong\mathbb{F}_4$ where $e$ is the
identity of ${\cal P}$. In this case $C_\varphi$ is a (Hermitian)
self-dual code of length $c$ over the quaternary field ${\cal P}$
under the inner product $(u,v)=\sum_{i=1}^{c}{u_iv_i^2}.$ Since
the minimum distance of $E_\sigma(C)$ is at least 8, this
Hermitian  code should have minimum distance at least 4.

To classify the codes, we need additional conditions for
equivalence. That's why we use the following theorem:

\begin{theorem} {\rm \cite{Yo}}\label{thm:eq}
The following transformations preserve the decomposition and send
the code $C$  to an equivalent one:
\begin{itemize}
  \item[(i)] a permutation of the fixed coordinates;
  \item[(ii)] a permutation of the 3-cycles coordinates;
  \item[(iii)] a substitution $x\rightarrow x^2$ in $C_{\varphi}$ and
  \item[(iv)]  a cyclic shift to each 3-cycle independently.
\end{itemize}
   \end{theorem}

\subsection{Codes with an automorphism of type $3$-$(6,26)$}

The extremal self-dual $[44,22,8]$ codes having an automorphism of
type $3$-$(6,26)$ are considered in \cite{Stef-SO} but the author
didn't succeed to classify all codes. We do this classification
now. Generator matrices of the codes $C_\varphi$ and
$E_{\sigma}(C)^*$ are presented in \cite{Stef-SO}. In the same
paper, it is also proved that $C_\pi$ is a binary  self-dual
$[32,16,\ge 4]$ code with a generator matrix
\begin{equation*} G_{\pi} = \left(
\begin{array}{c c}
0&D\\ I_6& F\\
\end{array} \right)
\end{equation*}
    where $D$ generates a $[26,10,8]$
    self-orthogonal code ${\cal D}^*$, and ${\cal D} _F$ is its dual code. The code ${\cal D}^*$ cannot be self-complementary (see \cite{Stef-SO}).
    According to \cite{BBGO}, there are 1768 inequivalent $[26,10,8]$
    self-orthogonal codes. Using as $D$ generator matrices of those
    codes which are not self-complementary, we obtain the self-dual
    $[44,22,8]$ codes invariant under the given permutation. To test
    them for equivalence, we use the program \textsc{Q-Extension}
    \cite{Iliya-aut}. The weight enumerators of the constructed
    codes are listed in Table \ref{Table:c6}.

\begin{theorem}
There are exactly $15621$ self-dual $[44,22,8]$ codes having an
automorphism of type $3$-$(6,26)$.
\end{theorem}

\begin{table}[htb]
\centering \caption{Extremal self-dual $[44,22,8]$ codes having an
automorphism of type $3$-$(6,26)$} \vspace*{0.2in}
\label{Table:c6}
{\footnotesize
\begin{tabular}{|c|ccccccccccccc|}
\noalign{\hrule height1pt} 
$\beta$&14&15&16&17&18&19&20&21&22&23&24&25&26\\
\hline
$W_1$&-&-&-&-&4&16&33&31&59&62&82&79&72 \\
$W_2$&11&26&58&201&342&433&505&462&677&685&717&599& 611\\
\hline\hline
$\beta$&27&28&29&30&31&32&33&34&35&36&37&38&39\\
\hline
$W_1$&47&72&48&51&51&68&54&64&39&54&38&38&29\\
$W_2$& 463& 490& 452& 485&654&724&674&851&558&530&430&438&327\\
\hline\hline
$\beta$&40&41&42&43&44&45&46&47&48&49&50&51&52\\
\hline
$W_1$&32&32&28&35&66&49&51&41&40&33&39&29&33 \\
$W_2$&328&238&194&120&140&72&89&43&85&13&46&5&27\\
\hline\hline
$\beta$&53&54&55&56&57&58&59&60&61&62&63&64&65\\
\hline
$W_1$&17&24&8&18&4&15&4&7&1&5&1&2&3\\
$W_2$&5&21&6&11&-&15&1&6&1&7&-&2&-\\
\hline\hline
$\beta$&66&67&68&70&72&74&76&82&86&90&104&122&154\\
\hline
$W_1$&5& 2& 1& 2& 1&2&-&1& 1&1&-&1&-\\
$W_2$&1&-&1&1&3&4&2&2&1&1&1&-&1\\
\noalign{\hrule height1pt}
\end{tabular}
}
\end{table}

\subsection{Codes with an automorphism of type $3$-$(8,20)$}

Up to equivalence, a unique Hermitian quaternary $[8,4,4]$ code
 exists (see \cite{MWOSW}). So up to equivalence we have a unique subcode $E_{\sigma}(C)^*$.
The code $C_{\pi}$  is a binary self-dual $[28,14,\ge 4]$ code
with a generator matrix $G_{\pi}$ given in (\ref{matGpi})  where
$B$ and $D$ generate self-orthogonal $[8,k_1,\ge 4]$ and
$[20,k_1+6,\ge 8]$ codes, respectively. Since $0\le k_1\le 4$,
${\cal D}^*$ is a binary self-orthogonal $[20,6\le k_2\le 10,\ge
8]$ code. All optimal binary self-orthogonal codes of length 20
are classified in \cite{Stef-SO}. There are exactly 23
inequivalent $[20,6,8]$ self-orthogonal codes, four inequivalent
$[20,7,8]$ self-orthogonal codes, and a unique $[20,8,8]$
self-orthogonal code. Hence $k_1\leq 2$.

In the case $k_1=2$ we obtain only two inequivalent extremal codes
of length 44, both with weight enumerator $W_{44,2}$, respectively
for $\beta=68$ and $\beta=76$. For $k_1=1$, there exist 52
self-dual $[44,22,8]$ codes, and for $k_1=0$, the inequivalent codes
number 5399. Their weight enumerators are listed in Table
\ref{Table:c8}.

\begin{theorem}
There are exactly $5453$ self-dual $[44,22,8]$ codes having an
automorphism of type $3$-$(8,20)$.
\end{theorem}

\textbf{Remark:} The extremal self-dual $[44,22,8]$ codes
invariant under a permutation of type 3-(8,20) are considered
independently in \cite{Kim-c8}. The author of that paper
 has
classified all extremal self-dual codes which have an automorphism
of order 3 with 8 independent 3-cycles.

\begin{table}[htb]
\centering \caption{Extremal self-dual $[44,22,8]$ codes having an
automorphism of type $3$-$(8,20)$} \vspace*{0.2in}
\label{Table:c8}
{\footnotesize
\begin{tabular}{|c|ccccccccccccccc|}
\noalign{\hrule height1pt} 
$\beta$&8&9&10&11&12&13&14&15&16&17&18&19&20&21&22\\
\hline
$W_1$&-&-&-&-&-&-&2&-& -& 5&5&3& 9&16&8 \\
$W_2$&2&-& 3&10&8& 27&47&81& 157&174&330&395&442&481&560\\
\hline\hline
$\beta$&23&24&25&26&27&28&29&30&31&32&33&34&35&36&37\\
\hline
$W_1$&16&  28&16&69&36&39&27&60&29&55&26&34&15&25&15\\
$W_2$&442&432&307&298&140&172&79&69& 41&56&13&25&9&9 &6\\
\hline\hline
$\beta$&38&39&40&41&42&43&44&45&46&49&50&52&53&68&76\\
\hline
$W_1$&15&3&8& 4&6&1& 4& -&3&1&1&-&1&-&-\\
$W_2$& 18&-&9& 4&6&3&4& -& 3&-&3&1&-&1&1\\
\noalign{\hrule height1pt}
\end{tabular}
}
\end{table}

\subsection{Codes with an automorphism of type $3$-$(10,14)$
}\label{sect:c10}

In this case  $C_{\varphi}$ is a
Hermitian self-dual $[10,5,4]$ code and by \cite{MWOSW} is
equivalent to either $E_{10}$ or $B_{10}$.
  As in \cite{BYR}, we can fix the generator matrix of the subcode
 $E_{\sigma}(C)^*$ in the following two forms, respectively:
$${\footnotesize
\begin{pmatrix}
011011011011000000000000000000\cr
 101101101101000000000000000000\cr
 000000011011011011000000000000\cr
 000000101101101101000000000000\cr
 000000000000011011011011000000\cr
 000000000000101101101101000000\cr
 000000000000000000011011011011\cr
 000000000000000000101101101101\cr
 011000011000011000011000101110\cr
 101000101000101000101000110011\cr
 \end{pmatrix} \ {\rm and} \
\begin{pmatrix}
011011011011000000000000000000\cr
101101101101000000000000000000\cr
000011101110011000000000000000\cr
000101110011101000000000000000\cr
000000000000000011011011011000\cr
000000000000000101101101101000\cr
000000000000000000011101110011\cr
000000000000000000101110011101\cr
000011110101000000011110101000\cr
000101011110000000101011110000\cr
\end{pmatrix}}.$$

 The code $C_\pi$ has parameters $[24,12,\geq 4]$. There are exactly thirty inequivalent
 such codes,  namely $E_8^3$, $E_{16}\oplus E_8$, $F_{16}\oplus E_8$, $E_{12}^2$, and the indecomposable codes  denoted by
 $A_{24}, B_{24},\dots, Z_{24}$  in  \cite{CPS32}.
  All codes have minimum weight 4 except the extended Golay code $G_{24}$ with minimum weight 8 and the code $Z_{24}$ with minimum weight 6.
   We use the generator matrices of the codes given in
  \cite{PlessSl}.   For any weight 4 vector in $C_\pi$ at most two nonzero
  coordinates may be  fixed points. An examination of the vectors
  of weight 4 in the listed codes eliminates 23 of them. By investigation
  of all alternatives for a choice of the 3-cycle coordinates in
    the remaining codes  $G_{24}$, $R_{24}$, $U_{24}$, $W_{24}$, $X_{24}$, $Y_{24}$ and $Z_{24}$ we
  obtain, up to
  equivalence, all possibilities for the generator matrix of the code $C_\pi$.

  Let $C_\pi$ be $R_{24}$. There is a unique possibility for
  the choice of the 3-cycle coordinates up to equivalence.
   The generator matrix of $C_\pi$ in
  this case can be fixed in the form
  $$
 G_{\pi}(R_{24})={\footnotesize
 \begin{pmatrix}
 1100000000 \ 11000000000000\cr
0110000000 \ 01100000000000\cr 0001100000 \ 00011000000000\cr
0000110000 \ 00001100000000\cr 0000001100 \ 00000011000000\cr
0000000110 \ 00000001100000\cr 0000000011 \ 00000000110000\cr
1001000000 \ 10010000001111\cr 1000001000 \ 10000010001100\cr
0000001110 \ 00000000010110\cr 1110000000 \ 00000000000111\cr
0001110000 \ 00000000001110\cr
 \end{pmatrix}}.$$

Let $\tau$ be a permutation of the ten cycle coordinates in
$G_{\pi}(R_{24})$.
 Denote by $C^{\tau }$ the self-dual $[44,22]$ code determined by
 $C_{\varphi}$ and the matrix $ \tau (G_{\pi}(R_{24}))$.

 We consider the products of transformations (ii), (iii)
and (iv) from Theorem \ref{thm:eq} which preserve
 the quaternary code $C_{\varphi}$.
  Their permutation parts form a subgroup of the symmetric group $S_{10}$ which we
   denote by  $L$. Let  $S=Stab(R_{24})$ be the stabilizer of the automorphism group of
the code generated by $G_{\pi}(R_{24})$ on the set of the fixed
points. It is easy to prove that if  $\tau_1$ and $\tau_2$ are
permutations from the group
 $S_{10}$, the codes $C^{\tau_1}$ and $C^{\tau_2}$ are equivalent
 iff the double cosets $S\tau_1L$ and $S\tau_2L$
 coincide. In our case $Stab(R_{24})=\left<(7, 8)(9,10)\right.$, $(7,9,10)$,
$(7,9)(8,10)$, $(7,10)$, $(5,6)$, $(4,6,5)$, $(2,3)$, $(1,3,2)$,
$\left.(1,4)(2,5)(3,6)\right>$.

When $C_{\varphi}=B_{10}$ we found  in \cite{BYR} a subgroup of
the group $L$ generated by the permutations $(3,4)(8,9)$,
$(1,2)(3,4)$, $(1,3)(2,4)$, $(6,7)(8,9)$, $(6,8)(7,9)$ and
$(1,6)(2,7)(3,8)(4,9)(5,10)$. So we obtain four $[44,22,8]$ self-dual
codes: $C_{B_{10}}^{id }$, $C_{B_{10}}^{(567)}$,
$C_{B_{10}}^{(36754)}$ and $C_{B_{10}}^{(368574)}$. These codes
have weight enumerator $W_{44,1}$ with $\beta = 60, 33,$ 30 and 21
and automorphism groups of orders $2^7\cdot 3^4$, $2^4\cdot 3^3$,
72  and 48, respectively.

When $C_{\varphi}=E_{10}$ the group
$L=\left<(1,3,5,7,9)(2,4,6,8,10), (1,2)(3,4)\right.$, $\left.
(1,3)(2,4), (9,10)\right>.$ We obtain seven $[44,22,8]$ self-dual
codes $C_{E_{10}}^{\tau }$ for $\tau \in \{$id$,  (4,5,6,7),
$(4,5,7)(6,9,8), $(2,3,5,4)$, $(2,3,5,4)(6,7)$,
$(2,3,5,7,4)(6,9,8)$, $(6,7)\}.$
 These codes
have also weight enumerator $W_{44,1}$ with $\beta = 42, 30,$ 36,
24, 42, 30 and 21 and automorphism groups of orders
$2^{10}\cdot3$, 24, 192, 36, $2^7\cdot3^2$, again 24 and 720,
respectively.

In this way from all the
 cases for $C_\pi$ we constructed
 1865 inequivalent
 $[44,22,8]$ self-dual
 codes with weight enumerator $W_{44,1}$ for $\beta = 10,\dots,52$, 54, 55, 60, 62, 65 and 6873 codes  with weight enumerator $W_{44,2}$ for
 $\beta = 3,\dots,36$, 38, 42, 45, 46, 50 and 52. The
 calculations for these results were done with the GAP Version 4r4
 software system and the program \textsc{Q-Extension} \cite{Iliya-aut}.
The results are summarized in Tables \ref{Table:c10a} and
\ref{Table:c10}.

\begin{table}[htb]
\centering \caption{Extremal self-dual $[44,22,8]$ codes having an
automorphism of type $3$-$(10,14)$} \vspace*{0.2in}
\label{Table:c10a} {\small
\begin{tabular}{|c|c|c||c|c|c||c|c|c|}
\noalign{\hrule height1pt} 
  & $W_{44,1}$ & $W_{44,2}$&  & $W_{44,1}$ & $W_{44,2}$&   & $W_{44,1}$ & $W_{44,2}$ \\
   \hline
   $G_{24}, B_{10}$ & 3 & 12 & $U_{24}, E_{10}$ & 74 & 49 & $Y_{24}, B_{10}$ & 136 & 746 \\
   \hline
   $G_{24}, E_{10}$ & 6 & 25 & $W_{24}, B_{10}$ & 71 & 11 & $Y_{24}, E_{10}$ & 456 & 2764 \\
   \hline
   $R_{24}, B_{10}$ & 4 & - &  $W_{24}, E_{10}$ & 188 & 33 & $Z_{24}, B_{10}$ & 71 & 541 \\
   \hline
   $R_{24}, E_{10}$ & 7 & - &  $X_{24}, B_{10}$ & 161 & 224 & $Z_{24}, E_{10}$ & 207 & 1824 \\
   \hline
   $U_{24}, B_{10}$ & 29 & 19 & $X_{24}, E_{10}$ & 459 & 635 &  &  &  \\
\noalign{\hrule height1pt}
\end{tabular}
}
\end{table}

\begin{theorem}
There are exactly $8738$ inequivalent self-dual $[44,22,8]$ codes
having an automorphism of type $3$-$(10,14)$.
\end{theorem}

\begin{table}[htb]
\centering \caption{Extremal self-dual $[44,22,8]$ codes having an
automorphism of type $3$-$(10,14)$} \vspace*{0.2in}
\label{Table:c10} 
{\footnotesize
\begin{tabular}{|c|cccccccccccccc|}
\noalign{\hrule height1pt} 
$\beta$&3&4&5&6&7&8&9&10&11&12&13&14&15&16\\
\hline
$W_1$&-&-&-&-&-&-&-&1&2&11&49&63&25&114\\
$W_2$&1 &3&31&31&93&143&183&377&428&560&622&552&755&510\\
\hline\hline
$\beta$&17&18&19&20&21&22&23&24&25&26&27&28&29&30\\
\hline
$W_1$&97&51&159&134&71&157&99&63&129&81&49&90&61&41\\
$W_2$&411&585&270&223&321&145&96&176&35&  71&  64&32&13& 57\\
\hline\hline
$\beta$&31&32&33&34&35&36&37&38&39&40&41&42&43&44\\
\hline
$W_1$&55&31&28&41&21&16&22&21&11&14&11&10&12&4 \\
$W_2$& 7&23&16&11& 3&8&-&9& -& -&-&4&-&-\\
\hline\hline
$\beta$&45&46&47&48&49&50&51&52&54&55&60&62&65&\\
\hline
$W_1$ &4& 4&1&1& 1& 2& 1& 2&1&1&1&1&1&\\
$W_2$&1&1&-&-&-&1&-&1&-&-&-&-&-&\\
\noalign{\hrule height1pt}
\end{tabular}
}
\end{table}

\subsection{Codes with an automorphism of type $3$-$(12,8)$}

In this case  $C_{\varphi}$ is a quaternary Hermitian self-dual
code of length 12 with minimum weight at least 4. There exist
exactly five inequivalent quaternary self-dual $[12,6,4]$ codes,
denoted by $d_{12}$, $2d_6$, $3d_4$, $e_6\oplus e_6$, and
$e_7+e_5$ in \cite{MWOSW}.

The code $C_\pi$ is a  binary self-dual $[20, 10, \geq 4]$ code.
There are exactly seven such codes, namely $d_{12}+d_8$,
$d_{12}+e_8$, $d_{20}$, $d_4^5$, $d_6^3+f_2$, $d_8^2+d_4$, and
$e_7^2+d_6$ \cite{CPS32}. 
Each choice for the fixed points can lead to a different subcode
$F_\sigma(C)$. We have considered all possibilities for each of
these seven codes, and found exactly 7 inequivalent codes for
$d_{12}+d_8$, one code for $d_{12}+e_8$, one code for $d_{20}$, 10
codes for $d_4^5$, 26 codes for $d_6^3+f_2$, 18 codes for
$d_8^2+d_4$, and 3 codes for $e_7^2+d_6$. Denote these codes by
$H_{i,j}$, for $i=1, 2\dots, 7$.

By the method used in Section \ref{sect:c10}, considering the
permutation parts of the products of transformations (ii), (iii)
and (iv) from Theorem \ref{thm:eq} and the stabilizer of the
automorphism group of the codes $H_{i,j}$ on the fixed points, we
classified all codes up to equivalence. There are exactly 122787
inequivalent codes. Their weight enumerators are of type
$W_{44,1}$ for $\beta=$ 10, \dots, 68, 70, 72, 74, 82, 86, 90, 122
and of type $W_{44,2}$ for $\beta=$ 0, \dots, 56,  58, \dots, 62,
64, 66, 68, 70, 72, 74, 76, 82, 86, 90, 104, 154. The values
obtained for $\beta$ are listed in Table \ref{Table:c12}.
\begin{theorem}
There are exactly $122787$ inequivalent self-dual $[44,22,8]$
codes having an automorphism of type $3$-$(12,8)$.
\end{theorem}

\begin{table}[htb]
\centering \caption{Extremal self-dual $[44,22,8]$ codes having an
automorphism of type $3$-$(12,8)$} \vspace*{0.2in}
\label{Table:c12}
{\footnotesize
\begin{tabular}{|c|cccccccccccc|}
\noalign{\hrule height1pt} 
$\beta$&0&1&2&3&4&5&6&7&8&9&10&11\\
\hline
$W_1$&-&-&-&-&-&-&-&-&-&-&789&556\\
$W_2$&7&151&594&1434&2178&3468&5793&7034&6881&9434&10031&6906\\
\hline\hline
$\beta$&12&13&14&15&16&17&18&19&20&21&22&23\\
\hline
$W_1$&313&1915&1072&655&2141&1105&912&1770&1029&736&1338&666\\
$W_2$&8502&7975&5072&4805&5111&2549&2552&2438&1692&1176&1609&778\\
\hline\hline
$\beta$&24&25&26&27&28&29&30&31&32&33&34&35\\
\hline
$W_1$&642&731&511&382&568&286&286&263&236&161&179&99\\
$W_2$&773&745&532&311&484&204&242&169&217&65&176&32\\
\hline\hline
$\beta$&36&37&38&39&40&41&42&43&44&45&46&47\\
\hline
$W_1$&126&87&88&55&69&38&52&28&48&17&32&10\\
$W_2$&73&42&68&30&44&29&30&21&21&9&26&10\\
\hline\hline
$\beta$&48&49&50&51&52&53&54&55&56&57&58&59\\
\hline
$W_1$&18&7&19&9&15&5&7&3&11&4&9&5\\
$W_2$&14&7&17&3&15&4&9&6&13&-&11&1\\
\hline\hline
$\beta$&60&61&62&63&64&65&66&67&68&70&72&74\\
\hline
$W_1$&3&1&2&1&2&3&4&2&1&2&1&2\\
$W_2$&6&1&5&-&2&-&1&-&1&1&3&4\\
\hline\hline
$\beta$&76&82&86&90&104&122&154&&&&&\\
\hline
$W_1$&-&1&1&1&-&1&-&&&&&\\
$W_2$&2&2&1&1&1&-&1&&&&&\\
\noalign{\hrule height1pt}
\end{tabular}
}
\end{table}



\subsection{Codes with an automorphism of type $3$-$(14,2)$}

The code $C_\pi$ in this case is a self-dual $[16, 8, 4]$ code.
There are exactly three such codes, namely $d_8^{2}$, $d_{16}$,
and $e_8^2$ \cite{CPS32}. We consider their generator matrices in
the form \[ G_1=gen(d_8^{2})={\small \left(\begin{array}{@{}c@{}}
1000000011100000\\
0100000011010000\\
0010000000001110\\
0001000000001101\\
0000100011001011\\
0000010011000111\\
0000001010111100\\
0000000101111100\\
\end{array}\right)}, \ \ G_2=gen(d_{16})={\small \left(\begin{array}{@{}c@{}}
1111000000000000\\
0011110000000000\\
0000111100000000\\
0000001111000000\\
0000000011110000\\
0000000000111100\\
0000000000001111\\
0101010101010101
\end{array}\right)}, \] and
$G_3=gen(e_8^2)=\left(\begin{array}{@{}c@{}}
HO\\OH\end{array}\right)$, where $H=(I_4\vert J+I_4)$, $J$ is the
all-ones $4\times 4$ matrix
and $O$ is the $8\times8$ zero matrix.
We have to consider permutations on these generator matrices that
can lead to different subcodes $F_\sigma(C)$. From all
possibilities for each of these codes we have found exactly 7
different cases for $C_\pi$ which can produce inequivalent codes
$C$, namely $G_1$, $G_1^{(1,16)}$, $G_1^{(3,16)}$, $G_2$,
$G_2^{(1,16)}$, $G_3$, and $G_3^{(1,16)}$.

The code $C_\varphi$ is a quaternary Hermitian self-dual
$[14,7,4]$ code. There are exactly 10 such codes, namely $d_{14}$,
$2e_7$, $d_8+e_5+f_1$, $2e_5+d_4$, $d_8+d_6$, $2d_6+f_2$,
$d_6+2d_4$, $3d_4+f_2$, $2d_4+1_8$, and $q_{14}$ \cite{MWOSW}.

Again, considering the permutation parts of the products of
transformations (ii), (iii) and (iv) from Theorem \ref{thm:eq},
and the stabilizer of the automorphism group of the codes $C_\pi$
on the fixed points, we classified all codes up to equivalence.

When $C_\pi=d_{16}$  all codes have weight enumerators $W_{44,1}$
for $\beta$=11, 14, 17, 20, 23, 26, 29, 32, 35, 38, 41, 44, 53,
62, and 65. When $C_\pi=e_8^2$ the weight enumerators are
$W_{44,1}$ for $\beta$=10, 13, 16, 19, 22, 25, 28, 31, 34, 37, 40,
43, 46, 49, 52, and 58. Lastly, when $C_\pi=d_8^{2}$ we
constructed codes with weight enumerator $W _{44,2}$ for
$\beta=$1, 2, 4, 5, 7, 8, 10, 11, 13, 14, 16, 17, 19, 20, 22, 23,
25, 26, 28, 29, 31, 32, 34, 35, 37, 38, 40, 41, 43, 44, 46, 52,
and 55. The total number of all self-dual $[44, 22, 8]$ codes,
having an automorphism of type $3$-$(14,2)$ is 243927. The results
are presented in Tables \ref{Table:c14a} and \ref{Table:c14}.


\begin{table}[htb]
\centering \caption{Extremal self-dual $[44,22,8]$ codes having an
automorphism of type $3$-$(14,2)$} \vspace*{0.2in}
\label{Table:c14a} {\small
\begin{tabular}{|c|c|c|c|c|c|}
\noalign{\hrule height1pt} 
    &$d_{14}$&$2e_7$&$d_8+e_5+f_1$&$2e_5+d_4$&$d_8+d_6$\\
\hline
   $d_{16}$   &  7 &  33  & 66   & 26  & 144 \\
\hline
  $ e_8^2$    &  9 &  20  & 77   & 26  & 197 \\
\hline
   $d_8^{2+}$ & 114 & 876 & 2907 & 490 & 6148 \\
\hline \hline
              &$2d_6+f_2$&$d_6+2d_4$&$3d_4+f_2$&$2d_4+1_8$&$q_{14}$\\
\hline
   $d_{16}$   & 573   & 384 & 2040 & 1663 & 1191 \\
\hline
   $e_8^2$    & 735   & 496 & 2830 & 2225 & 1561 \\
\hline
   $d_8^{2+}$ & 25841 & 14639 & 84081 & 60246 & 34520 \\
\noalign{\hrule height1pt}
\end{tabular}
}
\end{table}

\begin{theorem}
There are exactly $243927$ inequivalent self-dual $[44,22,8]$
codes having an automorphism of type $3$-$(14,2)$.
\end{theorem}

\begin{table}[htb]
\centering \caption{Extremal self-dual $[44,22,8]$ codes having an
automorphism of type $3$-$(14,2)$} \vspace*{0.2in}
\label{Table:c14}
{\footnotesize
\begin{tabular}{|c|cccccccccc|}
\noalign{\hrule height1pt} 
$\beta$&1&2&4&5&7&8&10&11&13&14\\
\hline
$W_1$&-&-&-&-&-&-&704&984&1912&1537\\
$W_2$&4565&4374&21709&15796&35653&26242&33236&22914&21064&14322\\
\hline\hline
$\beta$&16&17&19&20&22&23&25&26&28&29\\
\hline
$W_1$&2006&1281&1447&1008&978&493&480&384&295&147\\
$W_2$&10879&6663&4407&3053&1866&992&621&521&344&152\\
\hline\hline
$\beta$&31&32&34&35&37&38&40&41&43&44\\
\hline
$W_1$&123&124&98&29&17&54&21&8&9&18\\
$W_2$&109&88&85&19&16&24&14&9&4&2\\
\hline\hline
$\beta$&46&49&52&53&55&58&62&65&&\\
\hline
$W_1$&7&2&3&1&-&1&1&2&&\\
$W_2$&4&-&4&-&2&-&-&-&&\\
\noalign{\hrule height1pt}
\end{tabular}
}
\end{table}

\subsection{All self-dual $[44,22,8]$ codes with an automorphism of order $3$}

Here we summarize all obtained results for the extremal self-dual codes of length 44 having an automorphism of order 3. To test the codes for equivalence, we used the program \textsc{Q-Extension}. The classification result is given in the following theorem.

\begin{theorem}
There are exactly $394916$ inequivalent self-dual $[44,22,8]$
codes having an automorphism of order $3$.
\end{theorem}

We list the number of the codes with different weight enumerators in
Table \ref{Table:p=3}. For $\beta\ge 67$, all codes have
simultaneously automorphisms of type 3-(12,8) and also automorphisms
of type 3-(6,26). This proves that the orders of the automorphism
groups of these codes are multiples of 9. We give these orders in
Table \ref{Table:beta68}. All codes with $\beta\ge 63$ have
automorphisms of type 3-(6,26). All seven codes with $\beta=0$ have
automorphisms of type 3-(12,8). The full automorphism group for four
of them is the cyclic group of order 3, and the other three codes
have automorphism groups of order 12.

\begin{table}[htb]
\centering \caption{All extremal self-dual $[44,22,8]$ codes having an
automorphism of order $3$} \vspace*{0.2in}
\label{Table:p=3}
{\footnotesize
\begin{tabular}{|c|ccccccccccc|}
\noalign{\hrule height1pt} 
$\beta$&0&1&2&3&4&5&6&7&8&9&10\\
\hline
$W_1$&-&-&-&-&-&-&-&-&-&-&1487\\
$W_2$&7&4713&4968&1435&23881&19271&5824&42768&33242&9617&43614\\
\hline\hline
$\beta$&11&12&13&14&15&16&17&18&19&20&21\\
\hline
$W_1$&1539&324&3860&2659&680&4248&2471&972&3385&2182&851\\
$W_2$&30231&9070&29668&19954&5666&16669&9965&3804&7898&5880&2440\\
\hline\hline
$\beta$&22&23&24&25&26&27&28&29&30&31&32\\
\hline
$W_1$&2523&1327&807&1428&1080&512&1051&558&431&515&504\\
$W_2$&4798&2963&2095&2250&1985&978&1465&870&849&965&1048\\
\hline\hline
$\beta$&33&34&35&36&37&38&39&40&41&42&43\\
\hline
$W_1$&266&396&201&213&176&205&95&131&86&92&79\\
$W_2$&761&1082&606&609&477&526&348&366&258&221&134\\
\hline\hline
$\beta$&44&45&46&47&48&49&50&51&52&53&54\\
\hline
$W_1$&115&66&86&47&55&39&55&34&43&18&28\\
$W_2$&151&73&102&44&87&15&51&5&30&5&23\\
\hline\hline
$\beta$&55&56&57&58&59&60&61&62&63&64&65\\
\hline
$W_1$&8&25&5&20&5&7&1&6&1&2&3\\
$W_2$&6&15&-&17&1&6&1&7&-&2&-\\
\hline\hline
$\beta$&66&67&68&70&72&74&76&82&86&90&104\\
\hline
$W_1$&5&2&1&2&1&2&-&1&1&1&-\\
$W_2$&1&-&1&1&3&4&2&2&1&1&1\\
\hline\hline
$\beta$&122&154&&&&&&&&&\\
\hline
$W_1$&1&-&&&&&&&&&\\
$W_2$&-&1&&&&&&&&&\\
\noalign{\hrule height1pt}
\end{tabular}
}
\end{table}



\section{Extremal Self-Dual Codes of Length 44 with an automorphism of
order 7}

If $\sigma$ is an automorphism  of a binary self-dual $[44,22,8]$
code of order 7, then $\sigma$ is of type $7$-$(3,23)$ or
$7$-$(6,2)$ \cite{FFTAHuff}.

Let $h_1(x)=(x^3+x+1)$ and $h_2(x)=(x^3+x^2+1)$. As
$x^7-1=(x-1)h_1(x)h_2(x)$, we have ${\cal P}=I_1\oplus I_2$, where
$I_j$ is an irreducible cyclic code of length 7 with parity-check
polynomial $h_j(x)$, $j=1,2$. According Lemma \ref{M1M2M3},
$C_{\varphi}=M_1\oplus M_2$, where $M_j=\{u\in C_{\varphi}\mid
u_i\in I_j, i=1,\dots,c\}$ is a linear code over the field $I_j$,
$j=1,2$, and $\dim_{I_1}M_1+\dim_{I_2}M_2=c$. The polynomials
$e_1=x^4+x^2+x+1$ and $e_2=x^6+x^5+x^3+1$ generate the ideals
$I_1$ and $I_2$ defined above. Any nonzero element of $I_j=\{0,
e_j, xe_j\dots, x^6e_j\}, j=1,2$ generates a binary cyclic
$[7,4,3]$ code. Since the minimum weight of the code $C$ is 8,
every vector of $C_{\varphi}$ must contain at least two nonzero
coordinates. Hence the minimum weight of $M_j$ is at least 2,
$j=1,2$.

The transformation $x\rightarrow x^{-1}$ interchanges $e_1$ and
$e_2$. The orthogonal condition (ii) from Theorem \ref{thm2}
implies that once chosen, $M_1$ determines $M_2$ and the whole
$C_{\varphi}$. So we can assume, without loss of generality, that
$\dim_{I_1}M_1\leq \dim_{I_2}M_2$, and we can examine only $M_1$.

\subsection{ Codes with an automorphism of type $7$-$(3,23)$}

Let $C$ be a binary self-dual $[44,22,8]$ code having an
automorphism of type $7$-$(3,23)$. Then we have
$\dim_{I_1}M_1+\dim_{I_2}M_2=3$. Since the minimum weight of $M_2$
is at least 2, we have $1\leq\dim_{I_1}M_1\leq\dim_{I_2}M_2\leq
2$. Hence $\dim_{I_1} M_1=1$ and $\dim_{I_2}M_2=2$. Then $M_2$ is
an MDS $[3,2,2]$ code over the field $I_2$ and according to
condition (ii) from Theorem \ref{thm2},
$M_1=\langle(e_1,e_1,e_1)\rangle$ and
$M_2=\langle(e_2,e_2,0),(0,e_2,e_2)\rangle$.

 In this case $C_\pi$ is a binary self-dual code of length 26. If $v=(1100\ldots0)\in C_\pi$ then
$\pi^{-1}(v)+(\phi^{-1}(e_2,e_2,0),00\ldots0)$ will be a codeword
from $C$ of weight 6 which contradicts the minimum weight of $C$.
Hence in the notations of Theorem \ref{structure}, $k_1=0, k_2=10, k_3=3$, and $gen~C_\pi=\left(%
\begin{array}{cc}
  0 & D \\
  E & F \\
\end{array}%
\right)$, where the matrix $D$ generates a $[23,10,\geq 8]$ binary
self-orthogonal code. There are three such codes and their
generator matrices are given in \cite{BBGO}.
We take $E=I_3,$ and we determine
 the matrix $F$ using the condition (iii) of Theorem
 \ref{structure}. For each of the three codes there is a unique possibility for
the matrix $F$, up to equivalence. We obtain the codes $C_{7,1}$
with weight enumerator $W_{44,1}$ for $\beta=122$, $C_{7,2}$ with
weight enumerator $W_{44,2}$ for $\beta=104$, and $C_{7,3}$ with
weight enumerator $W_{44,2}$ for $\beta=154$. The orders of their
automorphism groups are
$3251404800=2^{15}\cdot3^{4}\cdot5^{2}\cdot7^{2}$,
$116121600=2^{13}\cdot3^{4}\cdot5^{2}\cdot7$, and
$786839961600=2^{16}\cdot3^{4}\cdot5^{2}\cdot7^{2}\cdot11^{2}$,
respectively. 
All of these codes have automorphisms of order 5 and are known
from \cite{St42-44}.


\begin{theorem} There are exactly three inequivalent binary
$[44,22,8]$ codes having an automorphism of type $7$-$(3,23)$.
\end{theorem}

\subsection{Codes with an automorphism of type $7$-$(6,2)$.}

Let $C$ be a binary self-dual $[44,22,8]$ code having an
automorphism of type $7$-$(6,2)$. Now $C_\pi$ is a binary $[8,4]$
self-dual code equivalent either to $C_2^4$ or the extended
Hamming code $E_8$, generated by the matrices
$G_1=\left(I_4|I_4\right)$ and $G_2 =\left(I_4|J+I_4\right)$,
respectively  where $I_4$ is the $4\times 4$ identity matrix and
$J$ is the all-ones $4\times 4$ matrix.

In this case $\dim_{I_1}
 M_1+\dim_{I_2}M_2=6$ and $1\leq\dim_{I_1} M_1\leq\dim_{I_2} M_2\leq 5$. Hence $\dim_{I_1} M_1=1,2,$ or 3.

\textbf{Case I:} $\dim_{I_1}M_1=1, \dim_{I_2}M_2=5$. It follows
that $M_2$ is an MDS $[6,5,2]$ code, and
$M_1=\langle(e_1,e_1,e_1,e_1,e_1,e_1)\rangle$. If $C_\pi= C_2^4$,
then $C_\pi$ contains a codeword $v=(v_1,00)$ such that
$\mbox{wt}(v_1)=2$. Since $M_2$ is an MDS code, it contains a
codeword $w$ of weight 2 with the same support as $v_1$. But then
the codeword $\pi^{-1}(v)+(\phi^{-1}(w),00)\in C$ has weight 6 - a
contradiction. Therefore $C_\pi= E_8$. Fixing the codes $M_1$ and
$M_2$ and considering all binary codes equivalent to $E_8$, we
found only one $[44,22,8]$ code with weight enumerator $W_{44,1}$
for $\beta=38$ and $|Aut(C)|=8064$.

\textbf{Case II:} $\dim_{I_1}M_1=2, \dim_{I_2}M_2=4$. We can take
\[ gen(M_1)=\left(%
\begin{array}{cccccc}
e_1&0&\alpha_1&\alpha_2&\alpha_3&\alpha_4\\
0&e_1&\alpha_5&\beta_1&\beta_2&\beta_3\\
\end{array}%
\right),\] where $\alpha_i\in\{0,e_1\}, i=1,\cdots,5$, and
$\beta_i\in I_1, i=1,2,3$. Considering all such matrices we obtain
nine possibilities such that the minimum weight of $M_1$ is $\geq
2$, up to equivalence. Here $gen(M_1)$ is written in the form
$(I_2|A)$, where
$A$ is one of the following matrices:\\
\[
 A_1=\left(
  \begin{array}{cccc}
    e_1 & 0 & 0 & 0 \\
    e_1 & e_1 & e_1 & e_1 \\
  \end{array}
\right), \ \ A_4=\left(
  \begin{array}{cccc}
    e_1 & e_1 & 0 & 0 \\
    e_1 &  xe_1 & e_1 & e_1 \\
  \end{array}
\right), \ \ A_7=\left(
  \begin{array}{cccc}
    e_1 & e_1 & e_1 & 0 \\
    e_1 &  xe_1 & x^2e_1 & e_1 \\
  \end{array}
\right),\]
\[ A_2=\left(
  \begin{array}{cccc}
    e_1 & e_1 & 0 & 0 \\
    0 &  0 & e_1 & e_1 \\
  \end{array}
\right), \ \ A_5=\left(
  \begin{array}{cccc}
    e_1 & e_1 & e_1 & 0 \\
    0 &   e_1 & e_1 & e_1 \\
  \end{array}
\right), \ \ \ A_8=\left(
  \begin{array}{cccc}
    e_1 & e_1 & e_1 & 0 \\
    e_1 &  xe_1 & x^3e_1 & e_1 \\
  \end{array}
\right),\]
\[ A_3=\left(
  \begin{array}{cccc}
    e_1 & e_1 & 0 & 0 \\
    0 &  e_1 & e_1 & e_1 \\
  \end{array}
\right), \ \ A_6=\left(
  \begin{array}{cccc}
    e_1 & e_1 & e_1 & 0 \\
    0 &   e_1 & xe_1 & e_1 \\
  \end{array}
\right),  \ \ A_9=\left(
  \begin{array}{cccc}
    e_1 & e_1 & e_1 & e_1 \\
    e_1 &  xe_1 & x^2e_1 & x^3e_1 \\
  \end{array}
\right).\]

In the case $C_\pi= C_2^4$, denote by $A_i^{\tau}$ the
$[44,22,8]$ code determined by $(I_2|A_i)$ and $C_\pi=\tau(G_1)$.
There are 21 inequivalent codes, namely $A_1^{id}$,
$A_2^{(2,5,7,3,6)}$, $A_3^{(2,5,6)}$, $A_3^{(2,5,7,3,6)}$,
$A_3^{(2, 5, 4, 6)}$, $A_4^{(2, 5, 7, 3, 6)}$, $A_5^{id}$,
$A_6^{(2, 5, 6)}$, $A_6^{(4, 5, 6)}$, $A_6^{(4, 5, 7, 8, 6)}$,
$A_6^{(3, 5)(4, 6)}$, $A_7^{id}$, $A_7^{(2, 5, 7, 3, 6)}$,
$A_7^{(4, 5, 7, 8)}$, $A_7^{(3, 6, 4, 5, 7)}$, $A_8^{id}$,
$A_8^{(2, 5, 7, 3, 6)}$, $A_9^{id}$, $A_9^{(3, 6, 7)}$, $A_9^{(4,
6, 7, 8)}$, and $A_9^{(2, 5, 6)}$. The code $A_2^{(2,5,7,3,6)}$
has an automorphism group of order $786839961600$ and is
equivalent to the  code $C_{7,3}$ constructed above.

In the case $C_\pi= E_8$, denote by $B_i^{\tau}$ the $[44,22,8]$
code determined by $(I_2|A_i)$ and $C_\pi=\tau(G_2)$. There are 19
inequivalent codes, namely $B_1^{id}$, $B_2^{id}$, $B_3^{id}$,
$B_3^{(5,6)}$, $B_4^{id}$, $B_4^{(3, 7, 6, 8, 5)}$, $B_5^{id}$,
$B_6^{id}$, $B_6^{(6, 7, 8)}$, $B_6^{(5, 6)}$, $B_6^{(4, 5, 6, 7,
8)}$, $B_6^{(3, 7, 8, 6, 4, 5)}$, $B_7^{id}$, $B_7^{(6, 7, 8)}$,
$B_7^{(5, 6)}$, $B_8^{id}$, $B_9^{id}$, $B_9^{(5, 6)}$, and
$B_9^{(5, 6, 7)}$. The code $B_2^{id}$  is equivalent to
$C_{7,2}$, constructed in the previous section.

\textbf{Case III:} $\dim_{I_1}M_1=\dim_{I_2}M_2=3$. Then
\[gen(M_1)=\left(%
\begin{array}{cccccc}
e_1&0&0&\alpha_1&\alpha_2&\alpha_3\\
0&e_1&0&\alpha_4&\beta_1&\beta_2\\
0&0&e_1&\alpha_5&\beta_3&\beta_4\\
\end{array}%
\right),\] where $\alpha_i\in\{0,e_1\}, i=1,\cdots,5$, and
$\beta_i\in I_1, i=1,2,3,4$. There are 18 inequivalent codes $M_1$
with minimum weight $\geq 2$. We can fix the generator matrices
for $M_1$ and $M_2$ and consider all possibilities for $C_\pi$.

When $C_\pi= E_8$ we obtain 64 inequivalent codes with $W_{44,1}$
for $\beta=10$, 17, 24, 31, 38, 52, and 122. In the case $C_\pi=
C_2^4$ we obtain 87 inequivalent codes with $W_{44,2}$ for
$\beta=0$, 7, 14, 21, 28, 35, 42, 56, and 154. The codes with
$\beta=122$ and $\beta=154$ are equivalent to $C_{7,1}$ and
$C_{7,3}$, respectively.

\begin{table}[htb]
\centering \caption{Automorphism groups of self-dual $[44,22,8]$
codes for $C_{\pi}=E_8$ } \vspace*{0.2in} \label{Table:p7c6h}
{\footnotesize
\begin{tabular}{|c|cccccccccc|}
\noalign{\hrule height1pt} 
   $|Aut(C)|$&7&14&21&28&42&56&84&112&126&168\\
\hline
$\sharp$ codes&13&35&1&5&9&2&2&2&1&1\\
\hline\hline
$|Aut(C)|$&252&336&672&1344&2688&5040&5376&8064&64512&3251404800\\
\hline
$\sharp$ codes&1&2&2&1&1&1&1&1&1&1\\
\noalign{\hrule height1pt}
\end{tabular}
}
\end{table}

\begin{table}[htb]
\centering \caption{Automorphism groups of self-dual $[44,22,8]$
codes for $C_{\pi}=C_2^4$ } \vspace*{0.2in} \label{Table:p7c6c}
{\footnotesize
\begin{tabular}{|c|cccccccc|}
\noalign{\hrule height1pt} 
    $|Aut(C)|$&7&14&28&42&56&112&168&336\\
\hline
Number of codes&49&33&4&1&3&2&1&2\\
\hline\hline
$|Aut(C)|$&672&1344&2688&10752&21504&43008&786839961600& \\
\hline
Number of codes&1&2&3&1&1&3&1& \\
\noalign{\hrule height1pt}
\end{tabular}
}
\end{table}
\begin{table}[htb]
\centering \caption{Weight enumerators of self-dual $[44,22,8]$
codes having an automorphism of order $7$} \vspace*{0.2in}
\label{Table:p7all}
{\footnotesize
\begin{tabular}{|c|cccccccccc|}
\noalign{\hrule height1pt} 
$\beta$ in $W_{44,1}$&&&10&17&24&31&38&52&59&122\\
\hline
Number of codes&&&23&19&14&12&9&4&1&1\\
\hline\hline
$\beta$ in $W_{44,2}$&0&7&14&21&28&35&42&56&104&154\\
\hline
Number of codes&27&29&32&5&7&1&1&4&1&1\\
 \noalign{\hrule height1pt}
\end{tabular}
}
\end{table}

\begin{theorem} There are exactly $191$ inequivalent $[44,22,8]$ codes
having an automorphism of order $7$.
\end{theorem}

The orders of the automorphism groups of these codes are presented
in Tables \ref{Table:p7c6h} and \ref{Table:p7c6c}. The weight
enumerators of the constructed codes are listed in Table
\ref{Table:p7all}.

\section{Summary}

 The self-dual $[44,22,8]$ codes having
automorphisms of order 11 are classified in \cite{164} and
\cite{160}. The codes invariant under an automorphism of order 5
are presented in \cite{St42-44} and \cite{Stef-SO}. Summarizing
these classifications and the results from the previous sections,
we obtain the following theorem.

\begin{theorem}
There are exactly $395555$ inequivalent self-dual $[44,22,8]$ codes
having an automorphism of odd prime order.
\end{theorem}

All constructed codes with $\beta\ge 43$ have automorphisms of order
3. In Table \ref{Table:all} we list the number of codes having an
automorphism of odd prime order according to their weight
enumerators but only for these values of $\beta$ for which there are
also codes having automorphisms of order 5, 7 or 11, but not 3. For
the other values of $\beta$ the number of all extremal self-dual
codes having an automorphism of odd prime order is the same as in
Table \ref{Table:p=3}. 
We can send the generator matrices of the obtained codes by e-mail to everybody who is interested.

\begin{table}[htb]
\centering \caption{Self-dual $[44,22,8]$ codes having an
automorphism of odd prime order} \vspace*{0.2in}
\label{Table:all}
{\footnotesize
\begin{tabular}{|c|cccccccc|}
\noalign{\hrule height1pt} 
$\beta$&0&4&5&7&9&10&11&12\\
\hline
$W_1$&-&-&-&-&-&{\bf 1506}&1539&{\bf 397}\\
$W_2$&{\bf 54}&{\bf 23926}&{\bf 19293}&{\bf 42796}&{\bf 9658}&{\bf 43639}&{\bf 30237}&9070\\
\hline\hline
$\beta$&14&15&17&19&20&21&22&24\\
\hline
$W_1$&2659&680&{\bf 2549}&3385&2182&851&{\bf 2561}&{\bf 820}\\
$W_2$&{\bf 20026}&{\bf 5672}&9965&{\bf 7909}&{\bf 5888}&{\bf 2445}&{\bf 4802}&{\bf 2117}\\
\hline\hline
$\beta$&25&27&28&29&30&31&32&\\
\hline
$W_1$&1428&{\bf 528}&1051&558&431&{\bf 525}&{\bf 523}&\\
$W_2$&{\bf 2251}&978&{\bf 1470}&{\bf 872}&{\bf 852}&965&1048&\\
\hline\hline
$\beta$&34&35&37&38&42&44&&\\
\hline
$W_1$&396&201&{\bf 179}&{\bf 207}&{\bf 96}&115&&\\
$W_2$&{\bf 1090}&{\bf 607}&477&526&221&{\bf 153}&&\\
\noalign{\hrule height1pt}
\end{tabular}
}
\end{table}


Looking at the tables, one can notice that there is only one code
for $\beta=154$. This code has a large automorphism group - its
order is
$2^{16}\cdot3^{4}\cdot5^{2}\cdot7^{2}\cdot11^{2}=786839961600$. The
same is the situation with the codes for $\beta=122$ and
$\beta=104$. These two codes have automorphism groups of orders
$2^{15}\cdot3^{4}\cdot5^{2}\cdot7^{2}=3251404800$ and
$2^{13}\cdot3^{4}\cdot5^{2}\cdot7=116121600$, respectively. In Table
\ref{Table:beta68} we present the orders of the automorphism groups
of the codes with $\beta\geq 67$. All these orders are
multiples of $288=9\cdot 2^5$. Actually, all 12 codes with
automorphism groups of orders bigger than 400000 have weight
enumerators of both types with $\beta\ge 72$ and are given in Table
\ref{Table:beta68}. We list the number of
codes $C$ with full automorphism groups of orders
$6000<|Aut(C)|<400000$, $|Aut(C)|\neq 2^s$, in Table \ref{Table:paut}. The code with the largest
automorphism group (order 368640) which is not listed in Table
\ref{Table:beta68} has weight enumerator $W_{44,1}$ with $\beta=42$.
Actually, the full automorphism group for most of the codes (exactly 309666) is the cyclic group of order 3. These codes have weight enumerators of both types with $\beta\le 42$.

\begin{table}[htb]
\centering \caption{The orders of the automorphism groups of the self-dual $[44,22,8]$
codes with $\beta\geq 67$} \vspace*{0.2in}
\label{Table:beta68}
{\footnotesize
\begin{tabular}{|c|c|c|c|c|c|c|}
\noalign{\hrule height1pt} 
   $\beta$ &67&68&70&72&74&76\\
\hline
$\sharp$ codes&&&&&&\\
with $W_{44,1}$&2&1&2&1&2&-\\
\hline $|Aut(C)|$&2592&5184&13824&6912 &6912&-\\
&2304&&18432&&73728&\\
\hline
$\sharp$ codes&&&&&&\\
with $W_{44,2}$&-&1&1&3&4&2\\
\hline &&&&92160&69120&207360\\
$|Aut(C)|$&-&207360&69120&1105920&331776-2&165888\\
&&&&184320&14745600&\\
\hline\hline
$\beta$ &82&86&90&104&122&154\\
\hline
$\sharp$ codes&&&&&&\\
with $W_{44,1}$&1&1&1&-&1&-\\
\hline $|Aut(C)|$&7372800&1105920&2211840&-&3251404800&\\
\hline
$\sharp$ codes&&&&&&\\
with $W_{44,2}$&2&1&1&1&-&1\\
\hline
&663552&&&&&\\
 $|Aut(C)|$&2211840&1105920&14745600&116121600&-&786839961600\\
\noalign{\hrule height1pt}
\end{tabular}
}
\end{table}

\begin{table}[htb]
\centering \caption{Number of the self-dual $[44,22,8]$ codes with
$6000<|Aut(C)|<400000$} \vspace*{0.2in} \label{Table:paut} {\small
\begin{tabular}{|c|ccccccc|}
\noalign{\hrule height1pt} 
    $|Aut(C)|$&368640&331776&207360&184320&165888&98304&92160\\
\hline
Number of codes&1&2&2&1&1&1&1\\
\hline\hline
$|Aut(C)|$&73728&69120&64512&61440&55296&46080&43008 \\
\hline
Number of codes&2&2&1&2&1&1&3 \\
\hline\hline
$|Aut(C)|$&36864&34560&21504&18432&15552&13824&12288 \\
\hline
Number of codes&4&2&1&8&1&2&11 \\
\hline\hline
$|Aut(C)|$&11520&10752&10368&9216&8064&6912&6144 \\
\hline
Number of codes&1&1&1&6&1&6&35 \\
\noalign{\hrule height1pt}
\end{tabular}
}
\end{table}

Looking at the weight enumerators of the extremal codes of length 44
constructed up to now, the following open problems arise:

\begin{enumerate}
\item Prove that there are not extremal self-dual $[44,22,8]$ codes with weight
enumerator $W_{44,1}$ for $\beta=69$, $71$, $73$, $75,\dots,81$,
83, 84, 85, 87, 88, 89, $91,\dots,121$, or  $W_{44,2}$ for
$\beta=57$, $63$, $65$, $67$, $69$, $71$, $73$, $75$,
$77,\dots,81$, 83, 84, 85, 87, 88, 89, $91,\dots,103$,
$105,\dots,153$.

\item Are the constructed codes with weight enumerators
$W_{44,1}$ for $\beta=61$, 63, 68, 72, 82, 86, 90, 122, and
$W_{44,2}$ for $\beta=59$, 61, 66, 68, 70, 86, 90, 104, 154, the
unique examples for their weight enumerators?

\item Which of these codes have connections with combinatorial
designs?
\end{enumerate}

\section*{Acknowledgements}

The authors would like to acknowledge the many helpful
suggestions of the anonymous reviewers.
We also thank the Editor of this Journal.

The first author thanks the
Department of Algebra and Geometry at Magdeburg University
for its hospitality while this work was completed, and  the
Alexander von Humboldt Foundation for its support.


%

\end{document}